\documentclass[11pt]{amsart}
\usepackage{amsfonts}
\usepackage{amsmath}
\usepackage{amssymb}
\usepackage{color}
\usepackage{graphicx}
\usepackage{xspace}
 \font \eightrm=cmr8

 \newcommand{\nc}{\newcommand}

\nc{\surj}{\to\hskip -3mm \to}

\nc{\ignore}[1]{{}}
\nc{\mrm}[1]{{\rm #1}}
\nc{\dirlim}{\displaystyle{\lim_{\longrightarrow}}\,}
\nc{\invlim}{\displaystyle{\lim_{\longleftarrow}}\,}
\nc{\vep}{\varepsilon} \nc{\ep}{\epsilon}
\nc{\sigmat}{\widetilde\sigma}
\nc{\ostar}{\overline{*}}

\nc{\mchar}{\mrm{Char}}
\nc{\Hom}{\mrm{Hom}}
\nc{\id}{\mrm{id}}

\nc{\remark}{\noindent{\bf{Remark:}}}
\nc{\remarks}{\noindent{\bf{Remarks:}}}

 \nc{\delete}[1]{}
 \nc{\grad}[1]{^{({#1})}}
 \nc{\fil}[1]{_{#1}}

\nc{\BA}{{\mathbb A}} \nc{\CC}{{\mathbb C}} \nc{\DD}{{\mathbb D}}
\nc{\EE}{{\mathbb E}} \nc{\FF}{{\mathbb F}} \nc{\GG}{{\mathbb G}}
\nc{\HH}{{\mathbb H}} \nc{\LL}{{\mathbb L}} \nc{\NN}{{\mathbb N}}
\nc{\PP}{{\mathbb P}} \nc{\QQ}{{\mathbb Q}} \nc{\RR}{{\mathbb R}}
\nc{\TT}{{\mathbb T}} \nc{\VV}{{\mathbb V}} \nc{\ZZ}{{\mathbb Z}}
\nc{\Cal}[1]{{\mathcal {#1}}}
\nc{\mop}[1]{\mathop{\hbox {\rm #1} }\nolimits}
\nc{\smop}[1]{\mathop{\hbox {\eightrm #1} }\nolimits}
\nc{\mopl}[1]{\mathop{\hbox {\rm #1} }\limits}
\nc{\frakg}{{\frak g}}
\nc{\g}[1]{{\frak {#1}}}
\def \restr#1{\mathstrut_{\textstyle |}\raise-8pt\hbox{$\scriptstyle #1$}}
\def \srestr#1{\mathstrut_{\scriptstyle |}\hbox to
  -1.5pt{}\raise-4pt\hbox{$\scriptscriptstyle #1$}}
\nc{\wt}{\widetilde}
\nc{\wh}{\widehat}
\nc{\un}{\hbox{\bf 1}}
\nc{\redtext}[1]{\textcolor{red}{\tt #1}}
\nc{\bluetext}[1]{\textcolor{blue}{#1}}
\nc{\comment}[1]{[[{\tt {#1}}]] }
\nc{\R}{\mathbb R}
\nc\fleche[1]{\mathop{\hbox to #1 mm{\rightarrowfill}}\limits}
\def\semi{\mathrel{\times}\kern -.85pt\joinrel\mathrel{\raise
    1.4pt\hbox{${\scriptscriptstyle |}$}}}
\nc{\np}{/\hskip -2.3mm\pi}
\nc{\snp}{/\hskip -1.8mm\pi}

\def\ta1{{\scalebox{0.2}{ %%%%%%%%%%%%%%%%%%%%%%%%%%%%%%%%%\ta1
\begin{picture}(12,12)(38,-38)
\SetWidth{0.5} \SetColor{Black} \Vertex(45,-33){5.66}
\end{picture}}}}
\begin{document}

\title[Rota--Baxter Algebra]
      {Rota--Baxter Algebra\\[0.3cm] {\small{The Combinatorial Structure\\[0.1cm] of Integral Calculus}}}

\author{Kurusch Ebrahimi-Fard}
\address{ICMAT,
		C/Nicol\'as Cabrera, no.~13-15, 28049 Madrid, Spain.
		On leave from UHA, Mulhouse, France}
         \email{kurusch@icmat.es, kurusch.ebrahimi-fard@uha.fr}         
         \urladdr{www.icmat.es/kurusch}

\author{Fr\'ed\'eric Patras}
\address{Laboratoire J.-A.~Dieudonn\'e
         		UMR 7351, CNRS,
         		Parc Valrose,
         		06108 Nice Cedex 02, France.}
\email{patras@math.unice.fr}
\urladdr{www-math.unice.fr/$\sim$patras}	    

%%%%%%%%%%%%%%%%%%%%%%%%%%%%%%%%%%%%%%%%%%%%%%%%%%%%%%%%%%%%%%%%%%%
\date{April 3rd, 2013}
%%%%%%%%%%%%%%%%%%%%%%%%%%%%%%%%%%%%%%%%%%%%%%%%%%%%%%%%%%%%%%%%%%%

\begin{abstract}
Gian-Carlo Rota suggested in one of his last articles the problem of developing a theory around the notion of integration algebras, complementary to the already existing theory of differential algebras. This idea was mainly motivated by Rota's deep appreciation for Kuo-Tsai Chen's seminal work on iterated integrals. As a starting point for such a theory of integration algebras Rota proposed to consider a particular operator identity first introduced by the mathematician Glen Baxter. Later it  was coined Rota--Baxter identity. In this article we briefly recall basic properties of Rota--Baxter algebras, and present a concise review of recent work with a particular emphasis of noncommutative aspects. 
\end{abstract}

\maketitle


%%%%%%%%%%%%%%%%%%%%%%%%%%%%%%%%%%%%%%%%%%%%%%%%

\tableofcontents

%%%%%%%%%%%%%%%%%%%%%%%%%%%%%%%%%%%%%%%%%%%%%%%%

\section{Introduction}
\label{sect:intro}

During the period 1960--1972 P.~Cartier, G.C.~Rota and collaborators set out to develop a general theory of the algebraic and combinatorial structures underlying integral calculus. Nowadays it has been subsumed by what is better known as the theory of Rota--Baxter algebras. In spite of the fact that those algebras need not to be commutative, a large part of the principal results have been described in the context of commutative function algebras. Interest in the general, i.e., noncommutative, case, which includes among several other areas the integral and finite difference calculus for operator algebras, started around 2000. Originally motivated by problems from theoretical physics, more precisely, the renormalization program in quantum field theory (that is, the process by which divergent integrals can be made finite), the domain of applications of the noncommutative theory has grown steadily since then.

The reader will find in \cite{Guo,KuRoYa,Rota3} good surveys on the foundations of the theory of Rota--Baxter algebras. Neither do we strive for an exhaustive presentation of the field nor do we seek to include a complete list of references. The latter restriction is mainly due to the extensive developments that have taken place during the last ten years, which resulted in many ramifications. The present survey will instead try to offer insights on recent progress made in a precise direction --maybe the most significative one from the point of view of integral and differential calculus-- namely the extension to the noncommutative setting of the fundamental results of the commutative theory, such as the Bohnenblust--Spitzer identity or the construction of the standard Rota--Baxter algebra (a presentation due to Rota of the free Rota--Baxter algebra in terms of symmetric functions).

The article commences by recalling the physical as well as the probabilistic origins of the theory and presents the main classical examples of Rota--Baxter algebras that stem mostly (but not exclusively) from analysis. The classical Cartier--Rota (commutative) theory is then presented briefly. Links with symmetric functions, shuffle algebras and fundamental (Spitzer-type) identities are put forward. The theory is presented with a view toward its noncommutative extension; this approach is neither standard nor, in some sense, the most natural one in the commutative setting, but it has the advantage of making transparent the ideas underlying the transition from the commutative to the noncommutative realm.

We enter then the core of the subject, i.e., the theory of noncommutative Rota--Baxter algebras by mainly insisting on three ideas. First, the existence of noncommutative analogs of most of the fundamental commutative \it constructions\rm. This idea is illustrated by the close relationship existing between the structure of free Rota--Baxter algebras and noncommutative symmetric functions \cite{gelfand}. Then, the existence of noncommutative analogs of most of the fundamental commutative \it formulas\rm. Eventually, we would like to emphasize that calculus in Rota--Baxter algebras is ``generic'' for several theories, very much as tensor calculus is generic for the theory of iterated integrals and Lie algebras (in a sense that can be made rigorous as in \cite{Reutenauer}: technically, formulas in the tensor algebra can be shown to hold in arbitrary enveloping algebras of graded Lie algebras). Theories to which such a remark applies in the noncommutative framework include the one of iterated integrals and sums of operator valued functions (corresponding respectively to noncommutative shuffle and quasi-shuffle products) as well as of derivations and differential operators (corresponding to pre-Lie products).

Concretely, all this means that it is often convenient (at least according to our experience) to work within the setting of Rota--Baxter algebras when it comes to looking for universal formulas solving problems in integral calculus, difference calculus, dynamical systems, and so on.

\vspace{0.5cm}
{\bf{Acknowledgements}}: The first author is supported by a Ram\'on y Cajal research grant from the Spanish government. The work was supported by the project MTM2011-23050 of the Ministerio de Econom\'{i}a y Competitividad, Spain. \\

%%%%%%%%%%%%%%%%%%%%%%%%%%%%%%%%%%%%%%%%%%%%%%%%

\section{Rota--Baxter Algebras}
\label{sect:RBA}

\begin{quote} 
``Whereas algebraists have devoted a lot of attention to derivations, the algebraic theory of the indefinite integral has been strongly neglected. The shuffle identities are only the tip of an iceberg of algebra and combinatorics of the indefinite integral operator which remains unexplored''. (G.C.~Rota \cite{Rota4}) 
\end{quote}

%%%%%%%%%%%%%%%%%%%%%%%%%%%%%%%%%%%%%%%%%%%%%%%%

\subsection{Origin and Integral Calculus}
\label{ssect:Origin}

The Rota--Baxter identity and, together with it, the very notion of Rota--Baxter algebra first appeared in 1960 in the work of Baxter\footnote{Glen Baxter and not Rodney Baxter --the latter is known, among others for the Yang--Baxter equations.}~\cite{Baxter}. Baxter's paper originated from a result in probability theory due to Spitzer in 1956 \cite{Spitzer} on which we will comment later. Various proofs of the so-called Spitzer identity had been obtained before, but Baxter's approach succeeded to unveil in terms of the aforementioned identity the underlying algebraic and combinatorial structures.

About ten years later, Rota \cite{Rota1,Rota2,Rota3} succeeded to deduce Spitzer's identity using classical results from the theory of symmetric functions --indeed, he showed that Spitzer's identity is \it equivalent \rm to the Waring identity. Rota and his school advocated algebraic combinatorics, both mathematically and on epistemological grounds --an attitude that can be recognized in his introduction to \cite{Rota1}:

\begin{quote} 
``The spectacular results in the fluctuation theory of sums of independent random variables [...] have gradually led to the realization that the nature of the problem, as well as that of the methods of solution, is algebraic and combinatorial [...]. It is the present purpose to carry this algebraization to the limit: the result we present amounts to a solution of the word problem for Baxter algebras''.\\
\end{quote}

Another essential component of the mathematical foundation of the theory can be traced back to P.~Cartier \cite{Cartier}. Similarly to Rota, Cartier considered the word problem for Rota--Baxter algebras (the problem amounts to constructing explicitly a basis of the free Rota--Baxter algebra), but he suggested a rather different approach, that lead him to the notion of quasi-shuffle product. Much later this product were formalized by M.~Hoffman in \cite{Hoffman}. It underlies for example the combinatorial structure of multiple zeta values (MZVs) \cite{CartierMZV}.

\medskip

In itself, the definition of Rota--Baxter algebra is fairly elementary. Indeed, it just consists of an associative algebra $A$ (say unital, and from now on over a field $k$ of characteristic $0$), equipped with a linear operator $R: A \to A$ satisfying the Rota--Baxter relation:
\begin{equation}
\label{RBR}
	R(x)R(y)=R\big(R(x)y+xR(y)+\theta xy \big).
\end{equation}
The parameter $\theta$ is a fixed element of $k$, which we call the \it weight\rm . The sum of the three terms on which the map $R$ acts on the right hand side of the equation defines a new product an $A$. We denote it $\ast_\theta$ and call it \it double product \rm (of weight $\theta$):
\begin{equation}
\label{double}
	x \ast_\theta y := R(x)y + xR(y) + \theta xy.
\end{equation}
Observe that (\ref{RBR}) implies that $R(x \ast_\theta y )=R(x)R(y)$. One checks easily that this product is associative and that $R$ is still a Rota--Baxter operator for $\ast_\theta$ (that is, the Rota--Baxter relation still holds when one replaces the initial product on $A$ by the double product $\ast_\theta$). We will later see how useful this new product may be.

Let us first consider the case when $\theta=0$ and $A$ is commutative. The Rota--Baxter relation is then simply the usual integration by parts relation: take for example $A$ to be the algebra $C^\infty(\RR)$ of smooth functions on the real line and set $R(f)(t):=\int_0^t f(x)dx$. For two functions $f,g$ we obtain the weight zero Rota--Baxter identity:
\begin{equation}
\label{RBR0}
	R(f)R(g)=R\big(R(f)g+fR(g)\big).
\end{equation}
Even in this somewhat degenerate case the theory is not without interest since it embraces among others Chen's iterated integral calculus, Sch\"utzenberger's theory of shuffle algebras \cite{schutz} as well as a large part of Fliess's approach to control theory \cite{Fliess} (the one of chronological algebras in the sense of \cite{Kawski}). Indeed, from (\ref{RBR0}) follows in the weight zero case that:
\begin{eqnarray}
\label{shuffle}
	\lefteqn{R(f_1R(f_2 \cdots R(f_{n}) \cdots ))R(g_1R(g_2 \cdots R(g_m) \cdots ))} \nonumber\\
	&\hspace{6cm} =\sum R(h_1R(h_2 \cdots R(h_{n+m}) \cdots )),
\end{eqnarray}
where the sum runs over all the shuffles of the sequences $f_1, \ldots, f_n$ and $g_1, \ldots, g_m$, that is, those sequences $h_1, \ldots, h_{m+n}$ made of the $f_i$ and the $g_j$ where the partial order of the $f_i$ and the $g_j$ is preserved (for example, $x_1y_1x_2y_2y_3x_3x_4$ or $x_1x_2y_1y_2x_3x_4y_3$ are shuffles of the words $x_1x_2x_3x_4$ and $y_1y_2y_3$, but not $x_1x_4x_2y_1y_2x_3y_3$, since $x_4$ appears on the left of $x_2$ and $x_3$). When $R$ is the integral operator, we recognize in $R(f_1R(f_2 \ldots R(f_{n}) \ldots ))$ the iterated integral of functions $f_1, \dots ,f_n$, and (\ref{shuffle}) yields the well known --shuffle-- product of two such iterated integrals.

Equivalently, from a more abstract point of view, the (weight zero) double product $\ast_0=\ast$ splits into two components (known as the Eilenberg--MacLane and Sch\"utzenberger half-shuffles):
$$
	f \ast g = f\uparrow g + f\downarrow g
$$
where 
\begin{equation}
\label{halfSh}
	f\uparrow g:=fR(g),\ f\downarrow g:=R(f)g. 
\end{equation}
The products $\uparrow$ and $\downarrow$ satisfy the axiomatic characterization of shuffle algebras that appeared independently in algebraic topology and in Lie theory in the 1950's
 \cite{EM,schutz}:
\begin{equation}
\label{demishuffle0}
	a\downarrow b=b\uparrow a,\quad\   (a\uparrow b)\uparrow c=a\uparrow (b\uparrow c+c\uparrow b).
\end{equation}
The reader may check by herself that these relations are enough to insure that the product defined by the sum of compositions $\uparrow + \downarrow$ is associative as well as commutative \cite{schutz}.

These ideas are fairly universal and show up in many different contexts. Let $X$ be an alphabet and let us write $Sh(X)$ for the algebra freely generated by $X$, and the products $\uparrow$, $\downarrow$ modulo the relations (\ref{demishuffle0}); using the interpretation \ref{halfSh} of the two products $\uparrow$, $\downarrow$, $Sh(X)$ embeds as a subalgebra in the free Rota--Baxter algebra over $X$. In the theory of free Lie algebras as well as for MZVs, and actually in most application domains where the notion is relevant, $Sh(X)$ is called the free shuffle algebra over $X$ or the tensor algebra over $X$ equipped with the shuffle product \cite{Reutenauer}. In control theory it is (sometimes) called the free chronological algebra over $X$ (in the sense of Kawski), in the theory of operads, it is called the free Zinbiel algebra over $X$ \cite{lodV}, a tribute to Cuvier's dual notion of Leibniz algebras \cite{cuv1,cuv3}.

So far we have presented the simple weight zero case and related commutative algebraic structures. The extension to the non-zero weight case as well as the generalization to noncommutative algebras are somewhat more involved. In the following we will outline the general theory. Saying this, we should remark that for the sake of space we refrain from giving an exhaustive account. However, let us emphasize the underlying idea. Namely, Rota--Baxter algebra calculus is representative for a whole class of theories, including iterated integrals of scalar and operator valued functions --commutative and noncommutative shuffle products--, summation operations --commutative and noncommutative quasi-shuffle products--, differential operators and integration methods --pre-Lie products.

%%%%%%%%%%%%%%%%%%%%%%%%%%%%%%%%%%%%%%%%%%%%%%%%

\subsection{Spitzer's Identity}
\label{sect:Spitzer}

We consider now the case of non-zero weight $\theta \neq 0$. Note that if $R$ is a Rota--Baxter map of weight $\theta$, then the map $R':=\beta R$ for $\beta \in k$ is of weight $\beta \theta$. This permits to rescale the original weight  $\theta \neq 0$ to $\theta = +1$ (or $\theta = -1$). There exist several natural examples of zero and non-zero weight Rota--Baxter algebras:
\begin{itemize}

\item 
In classical fluctuation theory, where extrema of sequences of real valued random variables play a crucial role, the operator $R$ is defined on the characteristic function $F(t):=\EE[\exp({itX})]$ of a real valued random variable $X$:
$$
 	R(F)(t):=\EE[\exp({itX^+})],
$$
where $X^+:=\max(0,X)$. One can show that $R$ is a Rota--Baxter map of weight $\theta = -1$. Note that one may call it the ``historical'' example, which motivated Spitzers' and hence Baxter's original works \cite{Baxter,KuRoYa,Rota1,Spitzer}.

\smallskip

\item
Let $A$ be an associative algebra (not necessarily commutative) which decomposes into a direct sum of two (non-unital) subalgebras, $A=A^+\oplus A^-$. The two orthogonal projectors $\pi_\pm: A \to A^\pm$ are Rota--Baxter maps of weight $\theta=-1$. In fact, any projector which satisfies relation (\ref{RBR}) is of weight $\theta=-1$. The dominant example for such an algebra are the Laurent series $\CC[\epsilon^{-1},\epsilon]]$, which decompose into a ``divergent'' part $\epsilon^{-1}\CC[\epsilon^{-1}]$ and a ``regular'' part $\CC[[\epsilon]]$. Further below we will see its role in the modern approach to renormalization theory in perturbative quantum field theory. A nice example of a noncommutative Rota--Baxter algebra are $n \times n$ matrices together with the projector that maps a matrix $\alpha$ to the upper triangular matrix $\alpha^u$ defined by replacing all entries below the main diagonal by zeros.       

\smallskip

\item 
From the point of view of analysis and for the theory of discrete dynamical systems, the fundamental example of Rota--Baxter maps are summation operators. On functions $f$ defined on $\mathbb{N}$ and with values in an associative algebra, we define the summation operator $R(f)(n):=\sum_{k=0}^{n-1}f(k)$, which is a Rota--Baxter map of weight $\theta=1$. Note that it is the right inverse of the finite difference operator $\Delta(f)(n):=f(n+1) - f(n)$.       
 
\end{itemize}

\smallskip

Interest in Rota--Baxter algebras goes beyond its purely algebraic setting. Indeed, their relevance mainly stem from the possibility to describe interesting and universal combinatorial identities, which are useful, for instance, in the diverse set of examples ranging from integration and summation operators to projectors. One of canonical examples of such an universal combinatorial result is known as Spitzer's classical identity.

Let us return to the commutative case. Spitzer's classical formula allows to calculate the characteristic function of the extrema of a discrete process $S_i:=X_1+ \cdots +X_i$, defined as a sequence of partial sums of a sequence of independent and identically distributed real valued random variables $X_i$. In other words, we consider the new sequence of random variables $Y_i:=\sup(0,S_1,S_2, \ldots ,S_i)$. Spitzer's formula permits to calculate the characteristic function in terms of  the positive part of the partial sum, denoted $S_i^+$. Let $F$ be the characteristic function of $X$. In the setting of Rota--Baxter algebra we find:
$$
	\EE[\exp({itY_k})]=R(F\cdot R(F\cdot\  \cdots  R(F \cdot R(F)) \cdots ))=:R^{(k)}(F),
$$
and:
$$
	\EE[\exp({itS_k})]=R(F^k).
$$ 
From this we deduce an identity which is true in every commutative Rota--Baxter algebra of any weight~$\theta$, namely:
\begin{eqnarray}
\label{SpitzId}
	\log\Big(1+\sum\limits_{n>0}R^{(m)}(F)\Big)
					= R\Big(-\sum\limits_{n>0}\frac{(-1)^n\theta^{n-1}F^n}{n}\Big).
\end{eqnarray}

Let $\Omega_\theta'(F)$ denote the argument of the map $R$ on the right hand side of (\ref{SpitzId}). It is obvious that  $\Omega_\theta'(F)=\theta^{-1}\log(1 + \theta F)$. Moreover, the following identity can be derived: 
\begin{equation}
\label{magnus}
	\Omega_\theta'(F) = \frac{\ell_{\theta \Omega_\theta'(F)}}{\mathrm{e}^{\ell_{\theta \Omega_\theta'(F)}}-1}(F)
				    = F + \sum_{n > 0} \frac{B_n}{n!} \ell_{\theta \Omega_\theta'(F)}^{n}(F).
\end{equation}
Here $\ell$ is the left multiplication operator, $\ell_{x}(y):=xy$, and $B_n$ are the Bernoulli numbers, which appear in the series expansion of $\frac{x}{\mathrm{e}^x-1}$. Later we will see how this rewriting, which appears to be somewhat pointless in the commutative case, allows for a straightforward generalization of the above identity to noncommutative Rota--Baxter algebras. The key is Magnus' seminal work on the solution of linear differential equations \cite{Magnus} (which appeared in 1954, only two years before Spitzer's important paper). Magnus pioneered the problem of calculating the logarithm of the solution of a linear initial value problem written as a time-ordered exponential. We refer the reader to \cite{BCOR08} for a comprehensive review. 

Deeper into the combinatorial nature of commutative Rota--Baxter algebras goes the Bohnenblust--Spitzer formula, which is a multilinear identity similar to Spitzer's identity
\begin{equation}
\label{clBSp}
	\sum_{\sigma \in S_n}\!\! R\bigl(\cdots R(F_{\sigma(1)})F_{\sigma(2)} \cdots \bigr)F_{\sigma(n)} 
			= \sum_{\pi \in \mathcal{P}_n}\! \theta^{n -|\pi|}\! \sideset{}{^{*_{\theta}}}
							\prod_{\pi_i \in \pi}(b_i-1)! \, \Bigl(\prod_{j \in \pi_i}F_j\Bigr),   
\end{equation}
where $S_n$ is the symmetric group of order $n$. On the right hand side $\mathcal{P}_n$ denotes set partitions of $[n]$, and $|\pi|$ is the number of blocks of the partition $\pi \in \mathcal{P}_n$. The size of the $i$th block $\pi_i$ of $\pi$ is denoted $b_i$. Note that $\sideset{}{^{*_{\theta}}}\prod$ refers to the double product in (\ref{double}). The weight zero case, $\theta=0$, simplifies to: 
$$
	\sum_{\sigma \in S_n} R\bigl(\cdots R( R(F_{\sigma(1)})F_{\sigma(2)})\cdots \bigr) F_{\sigma(n)}
	= \sideset{}{ ^{*}}\prod^n_{i=1} F_i.
$$

With the goal to generalize identity (\ref{clBSp}) to noncommutative Rota--Baxter algebras, we again propose a slightly more involved rewriting. Recall the canonical cycle decomposition $c_1 \cdots c_k$ of a permutation $\sigma \in S_n$. Each cycle starts with its maximal element, and the cycles are written in increasing order of their first entries.  
For example, $(32)(541)(6)(87)$ is such canonical cycle decomposition. The $j$th cycle is denoted $c_j=(a_{j_0} a_{j_1}\ldots a_{j_{m_j-1}})$, where $m_j$ is the size of this cycle and $a_{j_0}>a_l$, ${j_1} \le l \le j_{m_j-1}$.    

Identity (\ref{clBSp}) rewrites:
\begin{equation}
\label{clBSpPerm}
	\sum_{\sigma \in S_n} R\bigl(\cdots R( R(F_{\sigma(1)})F_{\sigma(2)})\dots \bigr)F_{\sigma(n)} 
	=	\sum_{\sigma \in S_n}  \mathcal{D}^\theta_{\sigma}(F_1,\ldots,F_n), 
\end{equation}  
where for each permutation $\sigma$ written in its canonical cycle decomposition $c_1 \cdots c_k$ we define:
\begin{equation}
\label{lr}
	\mathcal{D}^\theta_{\sigma}(F_1,\ldots,F_n) 
	:= \sideset{}{^{*_{\theta}}}\prod_{j=1}^{k}\Big( \big(\prod_{i=1}^{m_j-1}r_{\theta F_{a_{j_i}}}\big)(F_{a_{j_0}}) \Big).
\end{equation}   
Now the operator $r$ denotes right multiplication $r_{x}(y):=yx$. Note that the second product on the right hand side is with respect to the composition of these multiplication operators. For instance, for the permutation $\sigma :=(43)(512) \in S_5$ we obtain:
$$
	\mathcal{D}^\theta_{\sigma}(F_1,\ldots,F_5) 
	= r_{\theta F_{3}}(F_{4}) *_\theta (r_{\theta F_{2}}r_{\theta F_{1}})(F_{5})  
	= \theta^3 (F_4 F_3)*_{\theta}(F_5 F_1F_2).
$$  
In the commutative setting, the product $(\prod_{i=1}^{m_j-1}r_{\theta F_{a_{j_i}}})(F_{a_{j_0}}) $ is independent of the order of the $a_{j_i},\ i\geq 1$. Therefore, several of the terms $\mathcal{D}^\theta_{\sigma}(F_1,\ldots,F_n)$ on the right hand side of (\ref{clBSpPerm}) may coincide. The resulting coefficients $\prod_{j=1}^{k}(m_j-1)!$ lead to the compact form (\ref{clBSp}) written in terms of set partitions. Indeed, choose a block $\pi_j$, say, of size $m_j$, of a partition $\pi$. Its first entry is fixed to be its maximal element. This block then corresponds to $(m_j-1)!$ different cycles of size $m_j$.

%%%%%%%%%%%%%%%%%%%%%%%%%%%%%%%%%%%%%%%%%%%%%%%%

\subsection{Combinatorial Approach: Cartier and Rota}
\label{ssect:Cartier-Rota}

The main idea underlying the works of Cartier and Rota is to deduce such Spitzer-type identities in the context of the free commutative Rota--Baxter algebra constructed over an alphabet $X$. Generally speaking, the universal property of free algebras then implies that such formulae are automatically valid in any, that is, in the category of such algebras. However, we will see that the two approaches chosen by Cartier and Rota are rather complementary.      

Cartier's construction \cite{Cartier} of a free commutative Rota--Baxter algebra is a forerunner of the modern notion of the quasi-shuffle algebra (see \cite{Hoffman}). In \cite{GK} Cartier's construction was made more explicit. Recall that the quasi-shuffle product is very present in the context of the theory of multiple zeta values (MZVs) and other generalizations of special functions (see \cite{CartierMZV}). 

Suppose we start with an alphabet $X$ equipped with a commutative monoid structure. We denote the monoid composition additively. Letters and words are elements of $X$, respectively sequences of letters. Shuffling of two words means to arrange their letters consecutively in such a way that the relative orders are preserved. Summing over all possible arrangements preserving the relative orders gives rise to the shuffle product.  For instance, the word $x_1x_4x_5x_2x_6x_3$ is part of the shuffle product of the two words $x_1x_2x_3$ and $x_4x_5x_6$. The quasi-shuffle product is based on the shuffle product, but it allows moreover for an additional operation using the underlying monoid structure of $X$ in the following way. In the monomials resulting from the shuffle product of two words any two consecutive letters may be summed provided that they are from different words. For example, the word $x_1x_4x_5x_2x_6x_3$ is included in the quasi-shuffle product of the words $x_1x_2x_3$and $x_4x_5x_6$ (i.e., no internal summation), but also the words $x_1x_4y_1y_2$ and $z_1x_5x_2y_2$, where $y_1:=x_5+x_2$, $y_2:=x_6+x_3$ and $z_1:=x_1+x_4$.

A concise mathematical description of the quasi-shuffle product of two words is given in terms of its recursive definition \cite{Hoffman}. Here we will present it in the context of Rota--Baxter algebra, since this makes the link between these two notions more transparent. Moreover, it indicates why the quasi-shuffle product appears to be naturally related to summation operators, such as for instance in the theory of MZVs, where the identity:
$$
	\zeta(p)\zeta(q)=\zeta(p,q)+\zeta(q,p)+\zeta(p+q)
$$
follows from the sum-representation:
$$
	\zeta(p):=\sum_{n>0}\frac{1}{n^p},\qquad\ \zeta(p,q):=\sum_{m>n>0}\frac{1}{m^pn^q}.
$$

The idea is similar to the weight zero case (i.e.~shuffles). Let us assume that $\theta=1$. The Rota--Baxter double product (\ref{double}) decomposes into three terms:
$$
	x \ast_1 y = R(x)y+xR(y)+xy.
$$
The first two terms correspond to the aforementioned ``half-shuffles''. The last part, $xy$, is just the original commutative and associative algebra product. Using the notations $x \uparrow y:=xR(y)$, $x\downarrow y:=R(x)y$ and $x \cdot y:=xy$ yields:
\begin{equation}
\label{demi-quasi-shuffle}
	x\downarrow y=y\uparrow x,\quad\   (x\uparrow y)\uparrow z=x\uparrow (y\uparrow z+z\uparrow y + y \cdot z).
\end{equation}
When compared to the ordinary shuffle product (weight zero), the last term, $y \cdot z$, explains the presence of additional terms in the quasi-shuffle product. This becomes already obvious at lower degrees when compared with the shuffle formula (\ref{shuffle}). In the following example we separate the parts coming from ``pure'' shuffling using brackets from the additional components that result from the weight $\theta=1$ terms.        
$$
	R(x)R(y)=[R(xR(y))+R(yR(x))]+[R(xy)], 
$$
\begin{eqnarray*}
	R(x)R(yR(z)) 	&=& [R(xR(yR(z)))+R(yR(xR(z)))+R(yR(zR(x)))]\\
				& & \qquad\ + [R(xyR(z))+R(yR(xz))].
\end{eqnarray*}
More general formulae for larger products follow from a recursive application of (\ref{RBR}). Note, however, that these identities are valid only in the commutative case.  

Rota's approach \cite{Rota1,Rota2} is rather different. It is based on the observation that for any commutative algebra $A$, the algebra $A^{\mathbf N}$ of functions from the non-negative integers into $A$, that is, sequences of elements of $A$ with pointwise product, is naturally equipped with a Rota--Baxter algebra structure defined by the operator:      
$$
	R(a_1,\ldots ,a_n,\ldots ) = (0,a_1,a_1+a_2,\ldots ,a_1+\cdots +a_n,\ldots ).
$$
Indeed, let $A$ be the algebra $k[x_1,\ldots ,x_n,\ldots ]$ of polynomials in an infinite number of indeterminates. One can then show that the Rota--Baxter subalgebra generated by the sequence $x:=(x_1,\ldots ,x_n,\ldots )$ is free. This establishes an immediate link to the well-understood classical theory of symmetric functions. In fact, the terms in Spitzer's classical identity, such as $R(x^n)$ and $R^{(n)}(x)=R(R^{(n-1)}(x)x)$, correspond in the theory of symmetric functions to power sums  ($\sum_ix_i^n$) and elementary symmetric functions ($\sum_{i_1<\cdots <i_n}x_{i_1} \cdots x_{i_n}$), respectively. From this one concludes that Spitzer's identity follows from Waring's formula, which expresses power sums as polynomials in elementary symmetric functions.

%%%%%%%%%%%%%%%%%%%%%%%%%%%%%%%%%%%%%%%%%%%%%%%%

\section{Theory of Noncommutative Rota--Baxter Algebras}
\label{sect:ncRBA}

Thanks to the particular way we presented the commutative case we do not have to dwell too much on the details of the theory of noncommutative Rota--Baxter algebras. Instead, we will try to sketch the techniques as well as theoretical ingredients, that permit a simple transition from the commutative to the noncommutative realm.   

In spite of multiple good reasons, theoretical as well as applied ones (e.g., integration and finite difference calculus with respect to operator algebra valued functions), interest in noncommutative Rota--Baxter algebra has been rather sporadic, compared to the amount of work that erupted shortly after Baxter's paper came out. In fact, only in the last 10 years or so the theory of noncommutative Rota--Baxter algebra has been developed systematically.          

Key in this renewed interest in the theory of noncommutative Rota--Baxter algebras is its link to renormalization in perturbative quantum field theory. Recall that the renormalization program permits to give a sense to integrals appearing in perturbative calculations in quantum field theory, which are otherwise divergent. The main idea is to extract and eliminate those parts that cause the divergencies in a coherent way, which is moreover compatible with the underlying physics \cite{CasKen,collins}. This extraction and elimination, i.e., the renormalization process, is combinatorial in nature. It was essentially described by Bogoliubov and others, and has been known in physics for a long time (in fact, Bogoliubov's and collaborators work \cite{BoPa} appeared around the same time as Magnus' and Spitzer's work). See \cite{Bogo} for one of the earliest textbooks on the subject. Renormalization theory plays a fundamental role in the standard model, and the spectacular results in particle physics that followed from it. Let us remark that from a mathematical point of view, despite the fact that other approaches to the renormalization problem, like, for instance, the Epstein--Glaser method \cite{EG} or Wilson's renormalization group \cite{Wilson}, allow for both a comprehensive understanding of the origin of divergencies as well as the very nature of the renormalization program, Bogoliubov's subtraction scheme (better known as BPHZ renormalization method \cite{collins}) still carries a veil of mystery.              

In the late 1990s, Kreimer, and Connes and Kreimer came up with a rather unexpected approach to the BPHZ renormalization method in terms of so-called combinatorial Hopf algebras \cite{CK1,Kreimer}. In their work on the recursive subtraction algorithm of Bogoliubov in the context of dimensional regularization together with minimal subtraction, they needed the subtraction scheme map to satisfy a particular set of identities; it was pointed out by C.~Brouder that these identities could be subsumed into a single relation, which was henceforth called {\it{multiplicativity constraint}} \cite{KreimerC}. A few years later it was realized that the identity proposed by Brouder coincided with the Rota--Baxter identity. Subsequent work of several people clarified the role of the Rota--Baxter relation in the aforementioned work of Connes and Kreimer. We refer the reader interested in details to Manchon's article in \it La Gazette des math\'ematiciens \rm\cite{Manchon1}, where he gives a concise and elegant account of the algebraic structures involved in the Connes--Kreimer theory, including the Rota--Baxter algebra aspect.      

For what considers this presentation, two aspects are of crucial importance. On the one hand, Bogoliubov's recursive algorithm and the Birkhoff decomposition derived from this by Connes and Kreimer are an instance of a general decomposition principle valid in any Rota--Baxter algebra, and known as Atkinson's decomposition. On the other hand, the work of Connes and Kreimer triggered considerable research into what is nowadays known as the theory of combinatorial Hopf algebras as well as into the theory of pre-Lie algebras. This lead also to a renewed interest in noncommutative Rota--Baxter algebras.

\smallskip

The rest of this work is organized as follows. In the following section we recall Atkinson's decomposition for Rota--Baxter algebras. The generalization of Spitzer-type formulas to noncommutative Rota--Baxter algebras is presented in Section \ref{ssect:ncSpitzer}. In this context, the notion of pre-Lie algebra is central. In Section \ref{ssect:ncCartier-Rota} we describe the noncommutative analogs of Cartier's and Rota's constructions. We indicate how for noncommutative Rota--Baxter algebras the theory of noncommutative symmetric functions plays a role analogue to that of symmetric functions for commutative Rota--Baxter algebras described by Rota. By using either examples or references to other works, we end this article with a brief outline of possible applications of Rota--Baxter algebras -- note that this last part mainly reflects the taste and interests of the authors.

%%%%%%%%%%%%%%%%%%%%%%%%%%%%%%%%%%%%%%%%%%%%%%%%

\subsection{Atkinson's Factorization}
\label{ssect:Atkinson}

In 1963, based on Baxter's and Spitzer's works \cite{Spitzer}, F.~Atkinson immersed himself into the factorization problem of Rota--Baxter algebras \cite{Atkinson}. However, he did not limit himself to commutative algebras. 

Let $A$ be a Rota--Baxter algebra of weight $\theta$. We denote by $\tilde{R}$ the operator $-\theta \id_A -R$. It is easy to show that it is a Rota--Baxter operator of weight $\theta$. Moreover, the identity $\tilde{R}(x \ast_\theta  y )=-\tilde{R}(x)\tilde{R}(y)$ holds. Note that for $\theta =0$, we have $\tilde{R}=-R$. For $x \in A$, Atkinson as well as Baxter were interested in the equations:
$$
	f=1+ \lambda R(fx), \qquad\  h = 1+ \lambda \tilde{R}(xh).
$$
The formal parameter $\lambda$ is introduced to circumvent any discussions regarding convergence or invertibilty issues. The solutions of the two equations are given in terms of the argument of the logarithm of the left hand side of the Spitzer identity (\ref{SpitzId}) for $f$ (here $F= \lambda x$), and an analogue series expansion for $h$. 

The identity $R(a)\tilde{R}(b) = R\big(a \tilde{R}(b)\big) + \tilde{R}\big(R(a)b\big)$ implies one of Atkinson's central results, i.e., the identity $fh=1-\lambda\theta fxh$ in $A[[\lambda]]$ which yields:
\begin{equation}\label{Atkins}
 1+\lambda \theta x = f^{-1}h^{-1}.
\end{equation}
Let us for a moment return to renormalization in the perturbative quantum field theory. We limit ourself to the algebraic perspective (without going into the details, for which we refer the reader to Manchon's paper \cite{Manchon1} and Collins' classical monograph \cite{collins}). The operator $R$ is called the subtraction scheme map and is defined analogous to the second example in Section~\ref{sect:Spitzer}: it is a projector into a subalgebra which, so to say, contains the divergences. The map $R$ is therefore an idempotent Rota--Baxter map of weight $\theta=-1$; it verifies moreover $R(1)=0$. We now omit the parameter $\lambda$. The term $1-x$ in Atkinson's identity (\ref{Atkins}) is the object to be renormalized, i.e., a series of divergent integrals.  The quantity $f$ encodes the divergences. Classically it is denoted $C$ and called counterterm in renormalization theory. The quantity $h^{-1}$ is the result of the renormalization process, and consists of a series of terms each of which is finite. Note that we neglect deliberately any 
aspects related to questions of summability, which is another problem, where simple algebra may not be enough to say much.

Bogoliubov's renormalization method is a recursive process, which permits to calculate $h^{-1}$ from $x$. From a Rota--Baxter algebra point of view, it is just an example of Atkinson's decomposition. From $f(1-x)=h^{-1}$ we deduce that:
$$	
	f=1+R(f)=1 + R(f-h^{-1})=1+R(fx),
$$
and:
$$
	h^{-1}=\tilde R(h^{-1})=\tilde R(f(1-x))=1-\tilde R(fx).
$$
In fact, $h^{-1}$ is in the image of the orthogonal projector $\tilde R$, and $f-1$ is in the image of $R$. Therefore, we have $R(h^{-1})=0$ and $\tilde R(f)=1$. Using the underlying series expansions (and, more precisely, the grading of the terms in these series), one shows that the two equations are coupled:
$$
	h^{-1}=1-\tilde{R}(fx),\quad\ f=1 + R(fx).
$$
Bogoliubov's algorithm then amounts to solving this system recursively degree by degree (here we indicate the degree of the terms in each series by an index $j$, such that for $j=1,2$ we find $h^{-1}_1=-\tilde{R}(x_1), f_1= R(x_1), h^{-1}_2=-\tilde{R}(f_1x_1+x_2), f_2=R(f_1x_1+x_2)$, etc.).

%%%%%%%%%%%%%%%%%%%%%%%%%%%%%%%%%%%%%%%%%%%%%%%%

\subsection{Noncommutative Spitzer Identity}
\label{ssect:ncSpitzer}

One of the natural questions is then to understand the expression on the left hand side of Spitzer's identity, i.e., the logarithm of the solution of the equation $f=1+ R(fx)$. 

Notice that, depending on the algebra $A$ and operator $R$ under consideration, one may either view this solution as the counterterm in renormalization theory, or, to mention two more familiar examples, the solution of a 1st order linear differential equation in a matrix algebra, or the linear fixed point equation of a discrete dynamical system. 
In particular, in view of the penultimate example, one would expect that a ``noncommutative Spitzer formula'' should solve, among others, the Baker--Campbell--Hausdorff problem (to compute the logarithm of the solution of a 1st order linear differential equation). We will see that this is indeed the case.

Key in the theory of noncommutative Rota--Baxter algebras is the existence of another algebraic structure beside the double product (\ref{double}). Indeed, any Rota--Baxter algebra comes equipped with an underlying pre-Lie (or Vinberg) algebra structure \cite{Cartier2,ChaLiv,Manchon2}. Recall that a (left) pre-Lie algebra is a vector space endowed with a (non-associative) product, which verifies an identity weaker than associativity, and known as (left) pre-Lie relation:
\begin{equation}
\label{pLidentity}
	(x\vartriangleright y)\vartriangleright z - x\vartriangleright (y\vartriangleright z)
	=(y\vartriangleright x)\vartriangleright z - y\vartriangleright (x\vartriangleright z).
\end{equation}
Right pre-Lie algebras are defined similarly. Pre-Lie algebras are Lie admissible, i.e., the bracket $[x,y]_\vartriangleright:=x\vartriangleright y - y \vartriangleright x$ satisfies the Jacobi identity: 
$$
	[x,[y,z]_\vartriangleright]_\vartriangleright+[y,[z,x]_\vartriangleright]_\vartriangleright
						+[z,[x,y]_\vartriangleright]_\vartriangleright=0.
$$ 
The canonical example of a pre-Lie algebra, which can be traced back to Caley's work on rooted trees, is given in terms of derivations. One verifies for example that the space of vector fields generated by the derivations $x^n\partial_x$ is a pre-Lie algebra with product $(x^n\partial_x)\vartriangleright(x^m\partial_x):=m \cdot x^{n+m-1}\partial_x$. 

Any Rota--Baxter algebra (of weight $\theta$) inherits automatically a (left) pre-Lie structure, which is defined similarly to the double product:
\begin{equation}
\label{pL}
	a\vartriangleright_\theta b:= R(a)b - bR(a) - \theta ba.
\end{equation}
In fact, both left and right pre-Lie products can be defined. We will focus on the former product, since $a\, {}_\theta\!\!\vartriangleleft\! b:=- b\!\vartriangleright_\theta\! a$ is right pre-Lie. Note that the Lie brackets associated to $\vartriangleright_\theta$ and $\ast_\theta$ identify: $[a,b]_{\vartriangleright_\theta}=[a,b]_{\ast_\theta}$. For $\theta = 0$, we have $a\!\vartriangleright_0\! b = [R(a),b]:=R(a)b-bR(a)$. Obviously, in the commutative case the pre-Lie product simplifies to the original product of the algebra (up to a scaling by the weight 
$\theta$).

To understand the role played by the pre-Lie product we return to the Bohnenblust--Spitzer formula, using the description presented in (\ref{clBSpPerm}), which appeared artificial in the commutative case. For $n=2$ we see quickly that:
$$
	F_1\ast_\theta F_2-\theta F_2F_1 = R(F_1)F_2  + R(F_2)F_1
$$
provided that the Rota--Baxter algebra is commutative. In the noncommutative case the slightly less obvious identity holds:  
\begin{eqnarray*}
	F_1 *_{\theta} F_2 + F_2 \vartriangleright_\theta  F_1 = R(F_1)F_2  + R(F_2)F_1.						
\end{eqnarray*} 

This simple procedure of replacing the algebra product by the pre-Lie product generalizes to all orders. Indeed, it is ``sufficient" to redefine in (\ref{clBSpPerm}) the operator $\mathcal{D}^\theta_\sigma$ by substituting in definition~(\ref{lr}) the operator {\makebox{$r_{\theta x}$ (right product by $\theta x$)} with the operator {\makebox{$r_{{\vartriangleright_\theta} x}$}}, again a right multiplication operator, but defined now in terms of the pre-Lie product, {\makebox{$r_{{\vartriangleright_\theta} x}(y):=y \vartriangleright_\theta x$}}. 

Using the example following the definition~(\ref{lr}), we substitute in the Bohnenblust--Spitzer formula the expression:
$$
	\mathcal{D}^\theta_{\sigma}(F_1,\ldots,F_5) = \theta^3 (F_4 F_3)*_{\theta}(F_5 F_1F_2),
$$  
by the one defined in terms of the pre-Lie product: 
\begin{eqnarray*}
	\mathcal{D}^{\vartriangleright_\theta} _{\sigma}(F_1,\ldots,F_5) 
	&=&r_{\vartriangleright_\theta  F_{3}}(F_{4}) *_\theta (r_{{\vartriangleright_\theta}  F_{2}}r_{\vartriangleright_\theta  F_{1}})(F_{5}) \\ 
	&=& (F_4 \vartriangleright_\theta F_3)*_{\theta}((F_5\vartriangleright_\theta  F_1)\vartriangleright_\theta F_2).
\end{eqnarray*}
With this redefinition, the Bohnenblust--Spitzer identity (\ref{clBSpPerm}) holds in any noncommutative Rota--Baxter algebra of weight $\theta$ \cite{EFMP}. The identity has even a $q$-analogue  \cite{NovThi}.

To approach the Baker--Campbell--Hausdorff problem generalized to Rota--Baxter algebra, it is useful to recall Spitzer's classical formula and its rewriting in terms of the generating series of the Bernoulli numbers (\ref{magnus}), and to compare it to Magnus classical solution of the problem. In his seminal 1954 paper \cite{Magnus}, Magnus considered the (say, matrix valued) solution $Y(t)$ of the linear differential equation $\dot{Y}(t)=A(t)Y(t)$, $Y(0)=1$.  He showed that the logarithm $\Omega(A)(t):=\log(Y(t))$ is the solution of the differential equation:
$$
	\dot{\Omega}(A)=\frac{ad_{\Omega(A)}}{\mathrm{e}^{ad_{\Omega(A)}}-1}(A),
$$
where $ad$ is the ordinary adjoint action ($ad_{x}(y)=[x,y]=xy-yx$). We remark that, since $\Omega(A)(0) = 0$, the fact that the indefinite Riemann integral is a weight zero Rota--Baxter map, implies that:
$$
	ad_{\Omega(A)}(A)=[\Omega(A),A]=  \dot{\Omega}(A)  \vartriangleright_0 A
				     = \ell_{\dot{\Omega}(A)\vartriangleright_0 }(A), 
$$
where $\ell_{x{\vartriangleright_0} }(y):=x\vartriangleright_0 y$. 

The generalization to arbitrary Rota--Baxter algebras unifies Spitzer's classical formula with Magnus' expansion thanks to the pre-Lie product. Let $A$ be a Rota--Baxter algebra of weight $\theta$. If $f$ is a solution of $f =1+\lambda R(f x)$ in $A[[\lambda]]$, then the element $\Omega'(\lambda x)$ defined such that $R(\Omega'(\lambda x))=\log(f)$ satisfies \cite{EFM2}:
\begin{equation}
\label{pLMag}
	\Omega'(\lambda x)
	=\frac{-\ell_{ \Omega'(\lambda x){\vartriangleright_\theta}}}{\mathrm{e}^{-\ell_{ \Omega'(\lambda x){\vartriangleright_\theta} }}-1}(\lambda x)
	=\lambda x+\sum_{n>0}\frac{(-1)^nB_n}{n!}\ell_{\Omega'(\lambda x){\vartriangleright_\theta}}^n(\lambda x),
\end{equation}
where $\ell_{x{\vartriangleright_\theta} }(y):=x\vartriangleright_\theta y$. This series is called pre-Lie Magnus expansion. Here, the prime notation shall remind the reader of Magnus' original differential equation. As a remark we mention that a similar approach applies to Fer's expansion \cite{EFM2}.  F.~Chapoton showed in \cite{Chap} that (\ref{pLMag}) is of significant interest in the context of the theory of Lie idempotents as well as the theory of operads. Let us make the first few terms of  $\Omega'(\lambda x)$ explicit:
\begin{eqnarray*}
 	\Omega'(\lambda x) &= & \lambda x + \frac 12 \lambda^2 x\rhd_\theta x 
	+\frac 14 \lambda^3  (x\rhd_\theta x) \rhd_\theta x
 	+\frac 1{12} \lambda^3 x \rhd_\theta  (x\rhd_\theta x) + \\
	&&\qquad - \lambda^4 \frac 18  ((x \rhd_\theta  x) \rhd_\theta x ) \rhd_\theta  x 
			- \lambda^4 \frac{1}{24}\Big( (x \rhd_\theta  (x \rhd_\theta  x )) \rhd_\theta  x\\
	&&\qquad \qquad + x \rhd_\theta  ((x \rhd_\theta  x ) \rhd_\theta x) +  (x \rhd_\theta x) \rhd_\theta  (x  \rhd_\theta x)  \Big) + \cdots.
\end{eqnarray*}
Observe that using (\ref{pLidentity}) the four order $\lambda^4$ terms reduce to: 
$$
	 -\frac{1}{6} \bigl((x \rhd_\theta x)  \rhd_\theta  x \bigr)  \rhd_\theta x
 - \frac{1}{12} x  \rhd_\theta \bigl((x  \rhd_\theta x)  \rhd_\theta  x \bigr).
$$

%%%%%%%%%%%%%%%%%%%%%%%%%%%%%%%%%%%%%%%%%%%%%%%%

\subsection{The word problem}
\label{ssect:ncCartier-Rota}

We have seen that the case of weight zero Rota--Baxter algebras is representative for the classical theory of integration (and differential) calculus. The pre-Lie product reduces to the ordinary Lie bracket twisted by the indefinite Riemann integral map,  $a\vartriangleright_0 b= R(a)b - bR(a)=:[R(a),b]$. Hence, Magnus classical expansion seen as the weight $\theta=0$ noncommutative analogue of Spitzer's formula clearly emphasizes the key role played by the pre-Lie structure in noncommutative Rota--Baxter algebras. Moreover, the simple generalization to the non-zero weight case allows for a coherent and natural passage to the full noncommutative analogue of Spitzer's formula.   

The classical {\it{word problem}} leads to another class of fundamental ideas and structures related Rota--Baxter algebras. We have seen that the noncommutative analogue of Cartier's theory is intimately related to noncommutative analogues of shuffle and quasi-shuffle algebras. A natural way to represent these structures is in terms of graphical, i.e., combinatorial objects, such as planar rooted trees \cite{AgMo,EFG1}.

On the other hand, the description of the noncommutative analogue of Rota's classical solution to the word problem \cite{EGBP,EFMP} is surprisingly straightforward. One simply replaces in Rota's original work the algebra of polynomials in an infinite number variables $k[x_1,\ldots,x_n,\ldots]$ by its associative analogue, i.e., the free associative algebra (or tensor algebra) over the alphabet  $X:=\{x_1,\ldots,x_n,\ldots\}$. Rota's results continue to hold in this noncommutative setting, that is, the operator $R(x_1,\ldots,x_n,\ldots):=(0,x_1,x_1+x_2,\ldots,x_1+\cdots+x_n,\ldots)$ still is a Rota--Baxter map, and the Rota--Baxter subalgebra generated by the sequence $x=(x_1,\ldots,x_n,\ldots)$ gives a presentation of the free Rota--Baxter algebra over one generator $x$.         

Beside the fact that this gives a simple answer to the word problem, it generalizes Rota's seminal insight, unveiling the link between free commutative Rota--Baxter algebras and symmetric functions, to one of the fundamental notions in modern algebraic combinatorics, i.e., the theory of noncommutative symmetric functions --the latter has been developed in the last 20 years by J.-Y.~Thibon and his collaborators (G.~Duchamp, F.~Hivert, J.-C.~Novelli \it inter alia\rm ) in a series of important articles from \cite{gelfand} to \cite{duc}.  

Indeed, recall that, if we calculate the terms on the left hand side of Spitzer's classical identity in the commutative case, i.e., $R^{(n)}(x)=R(R^{(n-1)}(x)x)$, using Rota's presentation of the free Rota--Baxter algebra on one generator, then we find a sequence of elementary symmetric functions $\sum_{0<i_1< \cdots <i_n<k}x_{i_1} \cdots x_{i_n}$ --this is the basis for Rota's proof of the Spitzer formula, as a corollary of Waring's identity. In the noncommutative case the same observation holds, but this time noncommutative variables enter the picture. As a result F.~Hivert's theory of quasi-symmetric functions in noncommutative variables applies. This approach shares the same advantages with Rota's original work. Indeed, it allows to apply a whole range of results and techniques from the theory of noncommutative symmetric functions (and variants of it, such as Solomon's descent algebra, which plays a crucial role in the theory of free Lie algebras \cite{Reutenauer}), see e.g. \cite{EGBP,EFMP}.

One may summarize this by the remark that this approach points in the direction of far reaching generalizations of Spitzer's formula beyond those which we have just sketched. This is of course not surprising, as it is coherent with the general observation that passage from the commutative to the noncommutative realm in general is non-unique. The theory of noncommutative symmetric functions \cite{gelfand} provides an illustrative example for this.

%%%%%%%%%%%%%%%%%%%%%%%%%%%%%%%%%%%%%%%%%%%%%%%%

\subsection{Applications and Perspectives}
\label{ssect:Appl}

We have seen that the theory of Rota--Baxter algebra naturally provides the setting to work abstractly with common notions such as indefinite integrals (or dually derivations), summation operators, splitting of algebras. The range of possible applications is therefore numerous.  We limit ourself to mention only two examples, one in the direction of differential equations, the other with respect to universal algebra.

%%%%%%%%%%%%%%%%%%%%%%%%%%%%%%%%%%%%%%%%%%%%%%%%

\subsubsection{Control Theory}
\label{sssect:Control}

Control theory of differential equations has been for a long time associated with the combinatorics of iterated integrals, and its formal variant, the combinatorics of words and shuffles. M.~Fliess' works \cite{Fliess2} play an eminent role in these developments at this interface between combinatorics and analysis.   

From the point of view of integration (recall that the indefinite Riemann integral is Rota--Baxter of weight zero), the underlying combinatorics is the one of weight zero Rota--Baxter algebras. In one dimension the commutative case can be used. Beyond that, i.e., in the case of operator valued functions, noncommutativitiy enters the picture.    

It is therefore not too surprising to discover certain ideas from noncommutative Rota--Baxter algebras in modern control theory. Indeed, chronological calculus is based on the notion of chronological algebra, which is the (weight zero) pre-Lie structure that derives from the Riemann integral in noncommutative algebras of matrix or operator valued functions. A similar remark applies also in the context of theory of numerical methods for differential equations.     

The work of A.~Agrachev and R.~Gamkrelidze serves as a remarkable example \cite{AG,Manchon2}. These authors systematically developed the theory of free chronological algebras, i.e., free pre-Lie algebras. This work allows for example to understand the relationship between the pre-Lie Magnus expansion and the discrete Baker--Campbell--Hausdorff formula (which allows to calculate $\log(\exp(x)\exp(y))$ when $x$ and $y$ do not commute:
\begin{equation*}
	BCH(x,y) := x+y +\frac 12[x,y] + \frac 1{12}([x,[x,y]]+[y,[y,x]]) + \cdots )
\end{equation*} 
in the context of Rota--Baxter algebra. Using their construction of the group of flows one can show, for a general Rota--Baxter algebra, that the product $l=fh$ of solutions $f$ and $h$ of the fixed point equations $f =1+ R(f x)$ and $h =1+\lambda R(h y)$, respectively, solves the equation $l =1+ R(l z)$, where:
$$
	z = x \bullet y := x +e^{-g_{\Omega'(x){\vartriangleright_\theta} }}y.
$$    
Here $\Omega'(x \bullet y) = BCH_{\ast_\theta}(\Omega'(x),\Omega'(y))$, and $BCH_{\ast_\theta}$ stands for the Baker--Campbell--Hausdorff formula defined in terms of the Lie bracket following from the Rota--Baxter double product $\ast_\theta$.

%%%%%%%%%%%%%%%%%%%%%%%%%%%%%%%%%%%%%%%%%%%%%%%%

\subsubsection{Algebra Splitting and Yang--Baxter Equations}
\label{sssect:SplitYB}

We have seen how in an associative algebra Rota--Baxter maps permit to define new algebraic structures. In the following we restrict ourself to the weight zero case. However, the reader should have no problem to generalize the results to arbitrary weight.   

Let $A$ be a weight zero Rota--Baxter algebra. The half-shuffles $x\uparrow y:=xR(y)$ and $x\downarrow y:=R(x)y$ combine to an abstract shuffle-like product on $A$. In the commutative case, one verifies that:
\begin{equation}
\label{demishuffle}
	a\downarrow b=b\uparrow a,\quad\   a\downarrow (b\downarrow c)=(a\downarrow b + b \downarrow a) \downarrow c.
\end{equation}

In the noncommutative case we do not longer have $a\downarrow b=b\uparrow a$, and the half-shuffles satisfy a system of axioms that characterizes noncommutative half-shuffles. In modern terms one says that $A$ is equipped with a dendriform algebra structure defined in terms of the compositions $\downarrow$ and $\uparrow$ \cite{Aguiar,Manchon1}:
\begin{eqnarray}
\label{demishuffleNC}
	(a\uparrow b)\uparrow c=a\uparrow (b\uparrow c+b\downarrow c)\\
	a\downarrow (b\uparrow c)=(a\downarrow b)\uparrow c\\
	a\downarrow (b\downarrow c)=(a\uparrow b+a\downarrow b)\downarrow c.
\end{eqnarray}

M.~Aguiar pointed out in \cite{Aguiar} that this splitting principle should not be limited to associative algebras --an idea that appears to be important from the point of view of universal algebra. In the same article, he established a link between weight zero Rota--Baxter operators and the equation: 
\begin{equation}
\label{ag}
	r_{13}r_{12}-r_{12}r_{23}+r_{23}r_{12}=0,
\end{equation}
known as the associative analog of the classical Yang--Baxter equation. He proved that if in an associative algebra $A$, the element $r=\sum u_i\otimes v_i \in A\otimes A$ satisfies (\ref{ag}), then the map $R_r(x):=\sum u_ixv_i$ is Rota--Baxter of weight zero (recall that $r_{12}$ is defined as $\sum u_i\otimes v_i\otimes 1 \in A\otimes A\otimes A$, and analogously for the other $r_{ij}$). 

Next we will show that there exists another link between Rota--Baxter operators and the Yang--Baxter equation. But before that, we look at the half-shuffle products $\uparrow$ and $\downarrow$. Let $V$ be a vector space equipped with a product $\ast$ (a bilinear map with values in $V$, for the time being we do not assume that it has particular properties). We can then define the notion of a (weight 0) Rota--Baxter algebra operator $R$ on $V$ for the $\ast$ product by relation:
$$
	R(x)\ast R(y)=R(R(x)\ast y+x\ast R(y)).
$$
When the product satisfies relations (e.g.~associativity, commutativity, Jacobi identity, pre-Lie or shuffle identities...), one would expect the ``half--products'' $x\uparrow y=x\ast R(y)$ and $x\downarrow y= R(x)\ast y$ to inherit remarkable properties from those of $\ast$. We know already that this is the case when $\ast$ is associative or commutative, but this is true for a wider class of theories. We refer to \cite{EFG,Bai,Uchino}, where these phenomena are studied in detail and limit here our attention to the particular case of Lie algebras, that fits nicely in the framework of the present article.

Let $L$ be a Lie algebra equipped with a Rota--Baxter map $R$ of weight zero. Hence, $R$ satisfies:
\begin{equation}
\label{ybc}
	[R(x),R(y)]=R([R(x),y]+[x,R(y)]),
\end{equation}
which is the --operator form of-- classical Yang--Baxter equation. This is enough for the bracket $[x,y]_R:=([R(x),y]+[x,R(y)])$ to be a Lie bracket \cite{STS}. Analogous to the associative case, this ``double'' Lie bracket $[-,- ]_R$ is the addition of the products $x\uparrow y:=[x,R(y)]$ and $x\downarrow y:=[R(x),y]$: $[x,y]_R=x\uparrow y+x\downarrow y$. However, the two products  $\uparrow$ and $\downarrow$ are right and left pre-Lie, respectively and the Lie bracket $[ -,-]_R$ can be written: 
$$
	[x,y]_R =x\uparrow y-y\uparrow x.
$$
In the case of a Rota--Baxter algebra $A$ of arbitrary weight, the corresponding Lie algebra $A_L$ is a Rota--Baxer Lie algebra for the map $R$. We recover in that case the same algebraic structures as described further above (that is, the pre-Lie structures arising from the associative and Lie Rota--Baxter structures identify).  

Let us return to the general case of an arbitrary weight $\theta$. We can go a bit further in the direction of Yang--Baxter equations and their link to Rota--Baxter algebras. The classical Yang--Baxter equation admits another variant, which plays an important role in the seminal work of Semenov-Tian-Shansky \cite{STS}. He introduced the modified Yang--Baxter equation:
\begin{equation}
\label{ybm}
	[B(x),B(y)]=B([B(x),y]+[x,B(y)]) - \theta^2 [x,y].
\end{equation}
It is sufficient to guarantee that the bracket $[x,y]_B:=\frac{1}{2}([B(x),y]+[x,B(y)])$ is a Lie bracket. One verifies quickly that in an  associative Rota--Baxter algebra $A$ of weight $\theta$, the associative analogue:
\begin{equation}
\label{modRBR}
	B(x)B(y)=B(B(x)y+xB(y)) - \theta^2 xy,
\end{equation}
is satisfied by the operator $B:=R - \tilde{R}$. As well, the operator $B$ satisfies in $A_L$ the modified Yang--Baxter equation. Using Semenov-Tian-Shansky's terminology, the map $B$ defines an associative double structure on $A$. We finally remark that the Rota--Baxter double product (\ref{double}) rewrites:  $x \ast_\theta y = \frac{1}{2}(B(x)y + xB(y))$.
 
\medskip
 
These various results and examples presented in this article do certainly not give an exhaustive picture of the  (existing or forthcoming) application domains of Rota--Baxter algebra techniques. However, we believe they render the flavor of the theory, and its power to attack many different problems, and contribute strongly to the development of ideas that are relevant both to the general theory of algebraic structures and to applications.

%%%%%%%%%%%%%%%%%%%%%%%%%%%%%%%%%%%%%%%%%%%%%
%%%%%%%%%%%%%%%%%%%%%%%%%%%%%%%%%%%%%%%%%%%%%
%%%%%%%%%%%%%%%%%%%%%%%%%%%%%%%%%%%%%%%%%%%%%

\end{document}